\newif\ifieee
\ieeefalse

\ifieee
    \documentclass[10pt]{ieeeconf}
    \IEEEoverridecommandlockouts
    \def\pp{}
\else
    \documentclass[12pt]{article}
    \def\pp{.}
\fi

\usepackage{graphicx,psfrag,amsmath,versions}
\usepackage{fullpage}
\newcommand{\BEAS}{\begin{eqnarray*}}
\newcommand{\EEAS}{\end{eqnarray*}}
\newcommand{\BEQ}{\begin{equation}}
\newcommand{\EEQ}{\end{equation}}
\newcommand{\BIT}{\begin{itemize}}
\newcommand{\EIT}{\end{itemize}}

\newcommand{\ie}{{\it i.e.}}

\newcounter{oursection}

\usepackage{array}
\usepackage{caption,subcaption}
\usepackage{circuitikz,tikz}
\usepackage{caption}

\title{
  Optimal Racing of an Energy-Limited Vehicle
}

\ifieee
    \author{
      Nicholas Moehle and Stephen Boyd
      \thanks{Nicholas Moehle is with the Mechanical Engineering Department, Stanford University.
                \texttt{nicholasmoehle@gmail.com}.}
\else
    \author{
      Nicholas Moehle%
      \thanks{Mechanical Engineering Department, Stanford University. \texttt{moehle@stanford.edu}}
    }
\fi

\bibliographystyle{alpha}

\begin{document}

\maketitle

\begin{abstract}
We consider the problem of controlling
a vehicle to arrive at a fixed destination
while minimizing a combination of energy consumption and travel time.
Our model includes vehicle speed and accelaration limits,
aerodynamic drag, rolling resistance, nonlinear engine losses,
and internal energy limits.
The na\"{i}ve problem formulation is not convex;
however, we show that a simple convex relaxation is tight.
We provide a numerical example,
and discuss extensions to vehicles with unconventional drivetrains,
such as hybrid vehicles and solar cars.
\end{abstract}

\section{Introduction}
In this paper
we consider how to control a vehicle's longitudinal dynamics
to arrive at a fixed destination
while minimizing the energy consumption and travel time.
For traveling long distances,
a good strategy is to maintain a constant speed.
In many cases, however, this is not practical.
Such cases might include a hybrid vehicle encountering obstacles,
such as traffic signals or stop-and-go traffic,
or a racing vehicle that must decelerate to make sharp turns.
For a solar car, changing speed might be \emph{desirable},
as the predicted availability of solar power changes.
%In these cases, the speed profile satisfy the dynamic constraints
%on the vehicle, while also taking into account energy usage.
In all these cases,
it is useful for the vehicle to quickly plan a dynamic speed profile,
that meets both the dynamic and energy requirements of the vehicle.
The goal of this paper is to show how to do this.

We first focus on the specific problem of 
minimizing the energy required to reach
a destination with a fixed travel time.
We then use this problem as a building block to
solve related problems,
such as the minimum-travel-time problem
(without regard to energy consumption),
and the minimum-energy problem
(without a fixed travel time).
We note that although our focus in this paper is on ground vehicles,
the same principles could certainly be applied to aircraft or ships.

Our approach is a combination of vehicle longitudinal control,
which involves planning the vehicle position and velocity,
and drivetrain control,
which involves planning the internal energy and power usage.
In our view, previous attempts to link these domains
have been hampered by the nonlinearity of the equation
$(1/2) m v^2 = K$ relating the velocity and the kinetic energy,
which makes it difficult to link the vehicle dynamics (related to $v$)
and the drivetrain energy dynamics (related to $K$).
One key observation is that the convex relaxation $(1/2) m v^2 \leq K$
is tight under the assumption that ``more velocity is better''.
Speed and acceleration limits
(which seem to contradict the ``more velocity is better'' assumption)
can still be handled indirectly as constraints on the kinetic energy instead of the velocity.

We conclude with a numerical example.
We show the Pareto trade-off curve between energy consumption 
and travel time;
our convex reformulation is the main tool for exploring this trade-off curve.
%This Pareto curve contains some interesting trajectories:
%the minimum-travel-time trajectory,
%and the minimum-energy trajectory (assuming no time limit).
We then mention some simple extensions,
including hybrid vehicles and solar cars.

\subsection{Previous work}

\paragraph{Vehicle drivetrain control\pp}
The past few years has seen an explosion of research
related to drivetrain control for hybrid vehicles,
and providing an overview is beyond the scope of this paper.
For a good review, see \cite{panday2014review}.
For an early approach based on convex optimization,
see \cite{tate2000finding}.
Fewer papers consider controlling the drivetrain
and the vehicle dynamics simultaneously;
see \cite{uebel2018optimal} and references therein.

\paragraph{Control along a fixed path\pp}
Minimum-time longitudinal vehicle control
is a special case of minimum-time trajectory generation 
over a fixed path.
A well-known convex formulation of this problem
is reviewed in \cite{lipp2014minimum}.
This technique is applicable to very general vehicle models,
and can include constraints on the speed and acceleration along
the path.
Unfortunately, these constraints must hold pointwise;
formulations involving integrals of the speed and acceleration
(such as those required to limit energy consumption)
result in nonconvex constraints in general.

\paragraph{Convex optimization\pp}
Convex optimization problems can be solved efficiently and reliably
using standard techniques \cite{boyd2004convex}.
Recently, much work has been devoted to solving
moderately-sized convex optimization problems
quickly (\ie, in milliseconds or microseconds), 
possibly on embedded platforms,
which enables convex-optimization-based control policies to be implemented 
at kilohertz rates 
\cite{wang2010fast,odonoghue2012splitting}.
In addition, recent advances in automatic code generation for convex optimization
\cite{mattingley2009automatic,chu2013code}
can significantly reduce the cost
and complexity of developing and verifying an embedded solver.

\section{Model}
We propose the following model of an energy-limited vehicle.
The vehicle operates over the time interval $[0,T]$, along a fixed path.

\paragraph{Vehicle dynamics\pp}
The position $x_t$ (measured along the path) 
and velocity $v_t$ of the vehicle at time $t$ are related by
\begin{align}
\dot x_t = v_t.
\label{e-dynamics}
\end{align}
The initial conditions are $x_0 = x^{\rm init}$ and $v_0 = v^{\rm init}$.

\paragraph{Speed constraints\pp}
The velocity has upper and lower limits, \ie,
\begin{align}
v^{\rm min}_t \leq v_t \leq v^{\rm max}_t.
\label{e-speed-limits}
\end{align}
These bounds may depend on time.
We assume that $v^{\rm min}_t \geq 0$,
\ie, the vehicle cannot move backward.

\paragraph{Acceleration limits\pp}
The acceleration at time $t$ cannot exceed $a^{\rm max}_t$:
\begin{align}
\dot v_t \leq a^{\rm max}_t.
\label{e-accel-limits}
\end{align}
This can be used to model tire traction limits.
These could change over time,
as the vehicle performs lateral maneuvers
or encounters varying road conditions.

\paragraph{Kinetic energy\pp}
The kinetic energy of the vehicle at time $t$ is $K_t$,
which is defined as
\begin{align}
K_t = \frac12 mv_t^2,
\label{e-energy}
\end{align}
where $m$ is the mass of the vehicle,
which is positive.
The kinetic energy changes according to
\begin{align}
\dot K_t = P^{\rm drv}_t - P^{\rm drag}_t -P_t^{\rm rr} - P^{\rm brk}_t,
\label{e-energy-dyn}
\end{align}
where $P^{\rm drv}_t$ is the power delivered by the vehicle drivetrain
and $P^{\rm brk}_t$ is the brake power, which must be nonnegative.
The power lost to drag is
\[
P^{\rm drag}_t = \frac12 \rho AC_D v_t^3.
\]
Here $\rho$ is the density of the air,
$A$ is the frontal area of the vehicle,
and $C_D$ is the drag coefficient,
all of which are positive.
The power lost to rolling resistance is
\[
P_t^{\rm rr} = C^{\rm rr} v_t^2,
\]
where $C^{\rm rr}$ is the (positive) rolling resistance constant.
(This corresponds to a rolling resistance force linear
in the vehicle velocity.)

\paragraph{Drivetrain\pp}
The drive power comes from an on-board energy source with internal energy $E_t$.
This value could represent the state of charge of a battery,
or the quantity of combustible fuel remaining.
(In the sequel we will refer to it as the battery energy.)
This value changes according to
\begin{align}
\dot E_t = -f^{\rm eng}(P^{\rm drv}_t),
\label{e-drivetrain}
\end{align}
where $f^{\rm eng}$ is the engine characteristic,
which encodes the motor efficiency at different operating points.
The domain of this function is the interval
$[P^{\rm drv,min}, P^{\rm drv,max}]$,
where $P^{\rm drv,min}$ and $P^{\rm drv,max}$
are the minimum and maximum drive powers.
We assume this function is increasing,
which encodes the fact that 
increasing drivetrain power requires increasing energy consumption.
We also assume it is convex,
which encodes decreasing incremental efficiency.
The battery energy has initial condition $E_0 = E^{\rm init}$.
It must also respect the energy limits
\begin{align}
E^{\rm min} \leq E_t \leq E^{\rm max},
\label{e-energy-limits}
\end{align}
where $E^{\rm min}$ and $E^{\rm max}$ are constants.

\section{Optimal control}
We would like to control the vehicle
to reach (or exceed) a desired position $x^{\rm end}$ by time $T$,
and to do so while minimizing the energy consumed.
This is formalized as an optimal control problem:
\begin{equation}
\begin{array}{ll}
\mbox{maximize} & E_T \\
\mbox{subject to} 
   & x_T \geq x^{\rm end} \\
   & \mbox{displacement dynamics (\ref{e-dynamics}),} \\
   & \mbox{speed limits (\ref{e-speed-limits}),} \\
   & \mbox{acceleration limits (\ref{e-accel-limits}),} \\
   & \mbox{kinetic energy definition (\ref{e-energy}),} \\
   & \mbox{kinetic energy dynamics (\ref{e-energy-dyn}),} \\
   & \mbox{internal energy dynamics (\ref{e-drivetrain}),} \\
   & \mbox{internal energy limits (\ref{e-energy-limits}).}
\end{array}
\label{e-opt-ctrl}
\end{equation}
The variables are
the displacement $x$,
the velocity $v$,
the kinetic energy $K$,
the drive power $P^{\rm drv}$,
the (nonnegative) brake power $P^{\rm brk}$,
and the internal energy $E$,
which are functions over the interval $[0,T]$.
The problem parameters include the following constants:
the terminal position $x^{\rm end}$,
the initial conditions $x_{\rm init}$, $v_{\rm init}$, and $E_{\rm init}$,
the mass $m$,
the air density $\rho$,
the frontal area $A$,
the coefficient of drag $C_D$,
the rolling resistance coefficient $C^{\rm rr}$,
and
the energy limits $E^{\rm min}$ and $E^{\rm max}$.
The problem parameters also include the
acceleration limit $a^{\rm max}$
and the speed limits $v^{\rm min}$ and $v^{\rm max}$,
all of which are scalar-valued functions defined on $[0,T]$,
as well as the engine characteristic $f^{\rm eng}$,
a scalar-valued function defined on $[P^{\rm drv,min}, P^{\rm drv,max}]$.

\subsection{Convex relaxation}
Although (\ref{e-opt-ctrl}) is not convex as stated,
a simple relaxation yields a convex problem;
in \S\ref{s-tightness}, we show that this relaxation is tight.

Before forming the relaxation,
we first reformulate some of the existing constraints.
(Note that the reformulations in this paragraph are lossless,
\ie, they are not relaxations of the original constraints.)
The drag power and the rolling resistance power can both be expressed 
in terms of the kinetic energy as
\[
P_t^{\rm drag} = (1/2)\rho AC_D (2K_t/m)^{3/2}, \qquad
\]
and
\[
P_t^{\rm rr} = 2C^{\rm rr} K_t/m,
\]
respectively.
We can then eliminate $P_t^{\rm drag}$, $P_t^{\rm rr}$, and $P_t^{\rm brk}$
in the energy dynamics equation (\ref{e-energy-dyn}) to obtain
\begin{align*}
\dot K_t \leq P^{\rm drv}_t - (1/2) \rho AC_D (2K_t/m)^{3/2} - 2C^{\rm rr} K_t/m.
\end{align*}
Similarly, the maximum speed constraint and
the acceleration limit can both be written
using the kinetic energy as 
\[
K_t \leq (1/2) m(v_t^{\rm max})^2
\]
and
\[
\dot K_t \leq \sqrt{2mK_t}a^{\rm max}_t,
\]
respectively.
(The latter follows from the fact that $\dot v_t = \dot K_t/\sqrt{2mK_t}$.)

With these reformulations in mind,
we begin relaxing some constraints of problem (\ref{e-opt-ctrl}).
We start by relaxing the energy definition (\ref{e-energy}) to inequality:
\begin{align}
K_t \geq \frac12 mv_t^2.
\label{e-relaxed-energy-def}
\end{align}
We keep the initial condition $K_0 = (1/2)m(v^{\rm init})^2$
as an equality constraint.
We also relax the internal energy dynamics to
\begin{align}
\dot E_t \leq -f^{\rm eng}(P^{\rm drv}_t).
\label{e-relaxed-energy-dyn}
\end{align}
and enforce the bounds 
\begin{align}
f^{\rm eng}(P^{\rm drv,min}) \leq -\dot E_t \leq f^{\rm eng}(P^{\rm drv,max}).
\label{e-power-bounds}
\end{align}
These bounds simply state that $-\dot E$ is in the range of $f^{\rm eng}$,
and are implied by (\ref{e-drivetrain}).
%As we show below, both of these relaxations are in fact lossless.

The relaxed problem is
\begin{equation}
\begin{array}{ll}
\mbox{maximize} & E_T \\
\mbox{subject to} 
   & x_T \geq x^{\rm end} \\
   & \dot x_t = v_t \\
   & x_0 = x^{\rm init} \\
   & K_0 = (1/2)m (v^{\rm init})^2 \\
   & v_t \geq v_t^{\rm min} \\
   & K_t \leq (1/2) m(v_t^{\rm max})^2 \\
   & \dot K_t \leq \sqrt{2mK_t}a^{\rm max}_t \\
   & K_t \geq (1/2) mv_t^2 \\
   \ifieee
      & \dot K_t \leq P^{\rm drv}_t - (1/2) \rho AC_D (2K_t/m)^{3/2} \\
           & \qquad\qquad - 2C^{\rm rr} K_t/m \\
   \else
      & \dot K_t \leq P^{\rm drv}_t - (1/2) \rho AC_D (2K_t/m)^{3/2}
            - 2C^{\rm rr} K_t/m \\
   \fi
   & \dot E_t \leq -f^{\rm eng}(P^{\rm drv}_t) \\
   & E_0 = E^{\rm init} \\
   & E^{\rm min} \leq E_t \leq E^{\rm max} \\
   %\ifieee
   %   %& f^{\rm eng}(P^{\rm drv,min}) \leq -\dot E_t \\
   %   %& -\dot E_t \leq f^{\rm eng}( P^{\rm drv,max}). \\
   %   %& -\dot E_t \leq f^{\rm eng}( P^{\rm drv,max}).
   %\else
   & f^{\rm eng}(P^{\rm drv,min}) \leq -\dot E_t \leq f^{\rm eng}( P^{\rm drv,max}).
   %\fi
\end{array}
\label{e-relaxed-opt-ctrl}
\end{equation}
The variables in some of the constraints are indexed by $t$;
these constraints must hold for all $t\in [0, T]$.
The decision variables are $x$, $v$, $K$, $P^{\rm drv}$, and $E$,
which are scalar-valued functions defined on $[0,T]$.

\subsection{Tightness of relaxation}
\label{s-tightness}
The relaxation given above is in fact tight,
in the following sense:
Given a solution $z = (x,v,K,P^{\rm drv},E)$
to (\ref{e-relaxed-opt-ctrl}),
a solution for (\ref{e-opt-ctrl}) is given by
$\tilde z = (\tilde x,\tilde v, K,\tilde P^{\rm drv},\tilde P^{\rm brk}, E)$,
which is defined as
\ifieee
    \begin{align*}
        \tilde x_t &= \int_0^t \tilde v_\tau \; d\tau, \\
        \tilde v_t &= \sqrt{2K_t/m}, \\
        \tilde P_t^{\rm drv} &= (f^{\rm eng})^{-1}(-\dot E_t), \\
        \tilde P_t^{\rm brk} &= \tilde P^{\rm drv}_t - (1/2) \rho AC_D (2K_t/m)^{3/2}
             \\ & \qquad\qquad\qquad - 2C^{\rm rr} K_t/m - \dot K_t.
    \end{align*}
\else
    \[
        \tilde x_t = \int_0^t \tilde v_\tau \; d\tau,
        \qquad
        \tilde v_t = \sqrt{2K_t/m},
        \qquad
        \tilde P_t^{\rm drv} = (f^{\rm eng})^{-1}(-\dot E_t),
    \]
    \[
        \tilde P_t^{\rm brk} = \tilde P^{\rm drv}_t - (1/2) \rho AC_D (2K_t/m)^{3/2}
                 - 2C^{\rm rr} K_t/m - \dot K_t.
    \]
\fi
Here, $(f^{\rm eng})^{-1}$ is the inverse
of $f^{\rm eng}$,
which exists because $f^{\rm eng}$ is increasing on its domain.
Note that $\tilde P_t^{\rm drv}$
is well defined, as $-\dot E_t$ is in the domain of the inverse function
(due to the constraint (\ref{e-power-bounds})).
Note that $\tilde v_t$ and $\tilde P_t^{\rm brk}$
are both well defined, as $K_t \geq (1/2) m v_t^2$ ensures
that $K_t$ is nonnegative.

\paragraph{Proof of tightness\pp}
Here we show that the new point
is optimal for (\ref{e-opt-ctrl}).
Recall that $z$ is optimal for 
(\ref{e-relaxed-opt-ctrl}), which is a relaxation of (\ref{e-opt-ctrl}).
Because $z$ and $\tilde z$ generate identical objective values for
(\ref{e-relaxed-opt-ctrl}) and (\ref{e-opt-ctrl}), respectively,
then if $\tilde z$ is in fact \emph{feasible} for (\ref{e-opt-ctrl}),
it is optimal as well.
We show feasibility of $\tilde z$ below.

%the optimal value of (\ref{e-relaxed-opt-ctrl}) is $E_T$.
%Because the original point is optimal for (\ref{e-relaxed-opt-ctrl}),
%and (\ref{e-relaxed-opt-ctrl}) is a relaxation of (\ref{e-opt-ctrl}),
%the optimal value of (\ref{e-opt-ctrl}) must be greater than or equal to $E_T$.
%Because our constructed point obtains objective value $E_T$ for (\ref{e-opt-ctrl}),
%if it is also feasible, then it is optimal.

Due to (\ref{e-relaxed-energy-def}),
we have $ \tilde v_t = \sqrt{2K_t/m} \geq v_t \geq v_t^{\rm min}$.
Because $\tilde v_t \geq v_t$,
we also have 
$\tilde x_T
= \int_0^T \tilde v_\tau \; d\tau 
\geq \int_0^T v_\tau \; d\tau 
= x_T 
\geq x^{\rm end}$.
We also have
$f^{\rm eng}(\tilde P_t^{\rm drv}) = -\dot E_t$, as required.
To verify that
$\dot K_t \leq \tilde P^{\rm drv}_t - (1/2) \rho AC_D (2K_t/m)^{3/2} - 2C^{\rm rr} K_t/m$
holds,
we need only show that $\tilde P_t^{\rm drv} \geq P_t^{\rm drv}$.
This follows from $f^{\rm eng}( P_t^{\rm drv} ) \leq -\dot E_t$.
In particular, because $f^{\rm eng}$ is increasing,
we can invert this relation to obtain
\begin{align*}
P_t^{\rm drv} 
&\leq (f^{\rm eng})^{-1}(-\dot E_t) \\
&= \tilde P_t^{\rm drv}.
\end{align*}

Nonnegativity of the brake power follows from
$\dot K_t 
\leq \tilde P^{\rm drv}_t - (1/2) \rho AC_D (2K_t/m)^{3/2} - 2C^{\rm rr} K_t/m$.
The other constraints,
such as the maximum velocity limit,
the acceleration limit,
the kinetic energy definition,
and the energy limits,
are true by definition of the constructed point $\tilde z$,
or follow trivially from feasibility of the original point $z$.

\section{Example}
\begin{figure*}
\begin{center}
\includegraphics[width=1\textwidth]{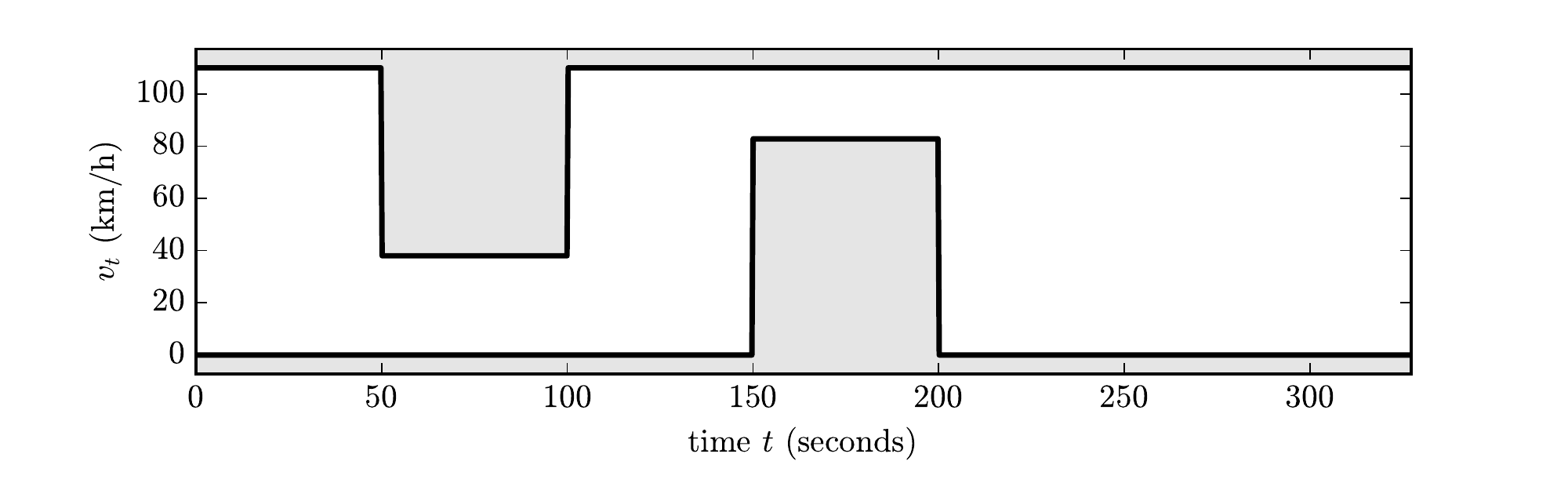}
\caption{The velocity limits $v_t^{\rm min}$ and $v_t^{\rm max}$ over time.
}
\label{f-limits}
\end{center}
\end{figure*}

We now present a simple numerical example.
The values of all scalar parameters are shown in table~\ref{t-params}.
The engine characteristic is
$f^{\rm eng}(p) = \alpha p^2 + \beta p + \gamma$
over the interval $[0,\infty)$.  
The parameters $\alpha$, $\beta$, and $\gamma$
are $0.005$ $\rm kW^{-1}$, $1$, and $5$ $\rm kW$, respectively.
The velocity limits, which are functions of time,
are shown in 
figure~\ref{f-limits}.
%The upper velocity limit $v^{\rm min}$ is around $110$ $\rm km/h$,
%except for a $50$ second period toward the beginning of the interval,
%during which it is reduced to around $40$ $\rm km/h$;
%the lower velocity limit $v^{\rm min}$ is zero,
%except for a $50$ second period toward the middle of the interval,
%when it is increased to around $80$ $\rm km/h$.
The maximum acceleration $a^{\rm max}$ was constant, 
with value $1$ $\rm m/s^2$. % TODO units
Several values of $T$ were considered,
to explore the trade-off between energy consumption and travel time.

\begin{table}
\begin{center}
\begin{tabular}{cc}
\hline
parameter & value \\
\hline
$m$ & $1500$ kg \\
$E^{\rm init}$ & $4000$ kJ \\
$E^{\rm min}$ & $0$ kJ \\
$E^{\rm max}$ & $4000$ kJ \\
$x^{\rm init}$ & $0$ m \\
$v^{\rm init}$ & $0$ m \\
$\rho$ & $1.22$ $\rm kg/m^3$ \\
$C_D$ & $0.35$ \\
$A$ & $2.3$ $\rm m^2$\\
$C^{\rm rr}$ & $0.005$ $\rm kN/(m/s)$\\
$x^{\rm end}$ & $5000$ m \\
\hline
\end{tabular}
\end{center}
\caption{
Parameter values.
}
\label{t-params}
\end{table}

To obtain a numerical solution,
problem (\ref{e-opt-ctrl}) was first discretized in time
by dividing the interval $[0,T]$ into $1001$ discrete points.
It was solved using \texttt{CVXPY} \cite{cvxpy},
with backend solver \texttt{ECOS} \cite{domahidi2013ecos}.

\begin{figure}
\begin{center}
\ifieee
    \includegraphics[width=1\columnwidth]{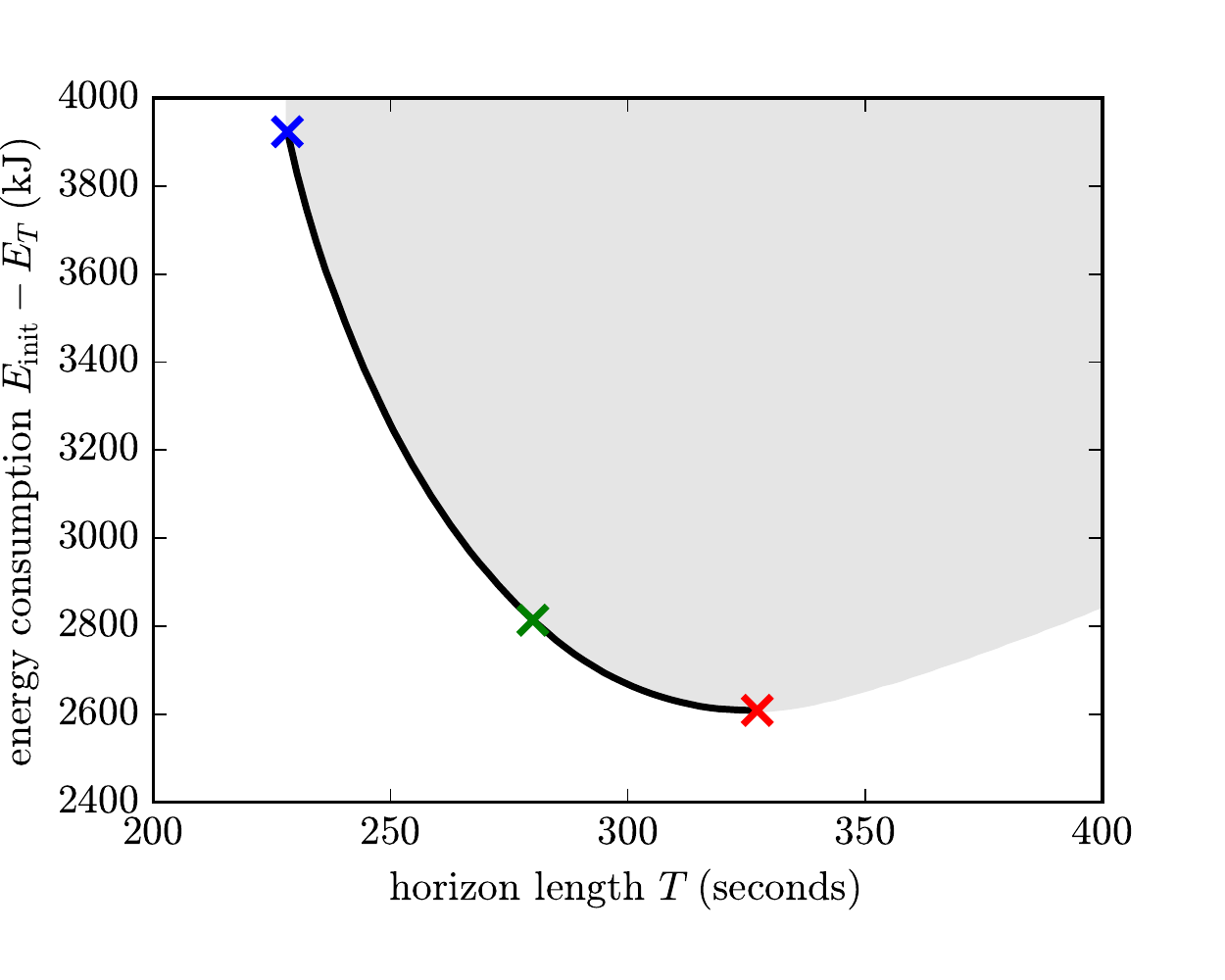}
\else
    \includegraphics[width=.7\columnwidth]{fig/pareto.pdf}
\fi
\caption{The Pareto optimal curve
trading off energy consumption and travel time.
The blue and red crosses
correspond to the minimum-time and minimum-energy trajectories,
and the green cross corresponds to a compromise between these two.
}
\label{f-pareto}
\end{center}
\end{figure}

\begin{figure*}
\begin{center}
\includegraphics[width=1\textwidth]{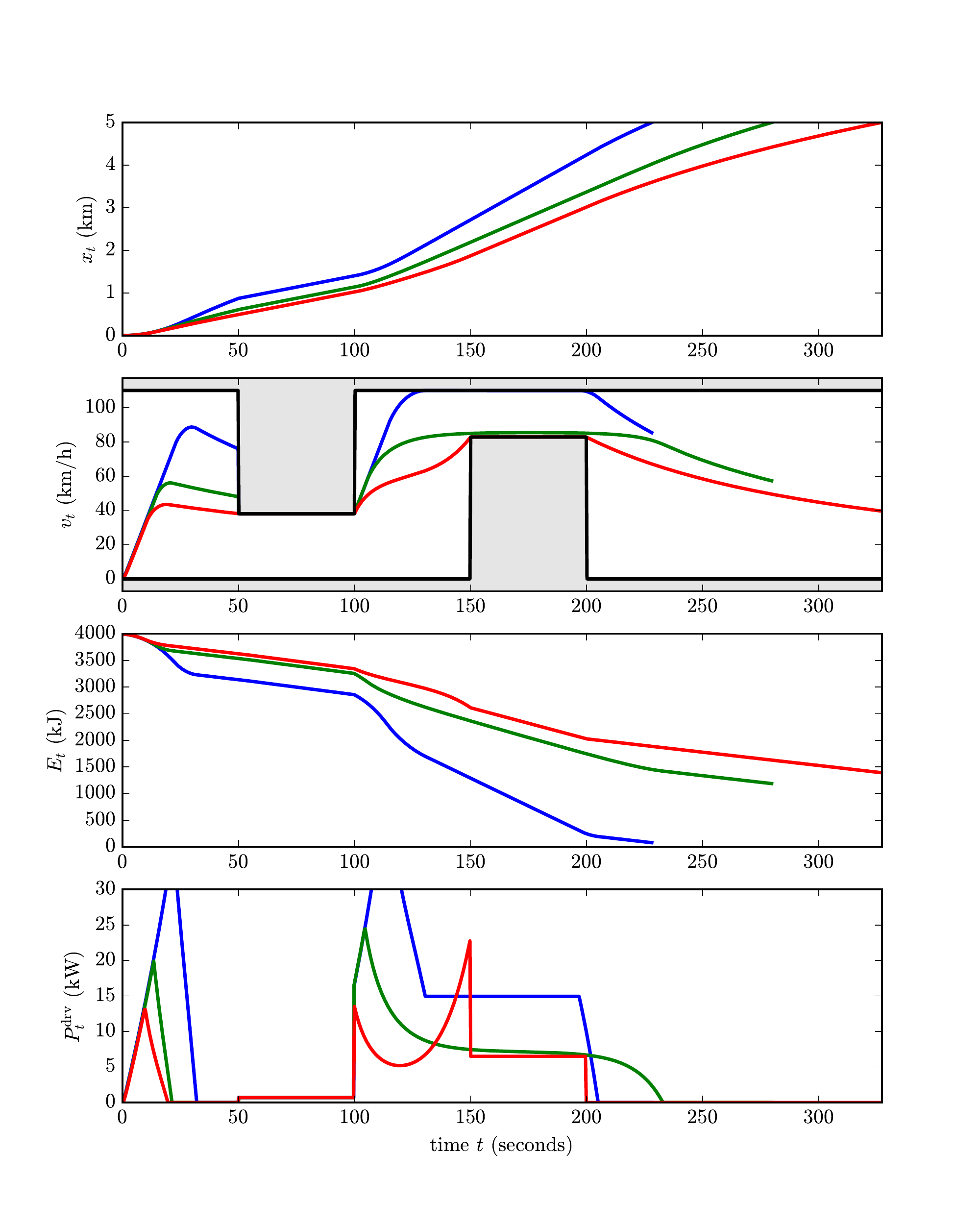}
\caption{Three optimal trajectories,
obtained by solving (\ref{e-opt-ctrl})
for different values of $T$.}
\label{f-trajectories}
\end{center}
\end{figure*}

\paragraph{Pareto curve\pp}
The trade-off between the total energy consumed
and the travel time is depicted in figure~\ref{f-pareto}.
The shaded region is
the set of possible pairs $(E^{\rm init} - E_T, T)$
of energy consumption and travel time
corresponding to feasible trajectories of (\ref{e-opt-ctrl}).
The Pareto curve
%which shows the trade-off between 
%minimizing energy consumption and travel time,
is shown as a bold line.
Trajectories correspending to any point on this line
can be computed by fixing $T$ and solving (\ref{e-opt-ctrl}).
The three colored crosses are specific Pareto-optimal trajectories,
which are shown in detail in figure~\ref{f-trajectories}.
%In particular,
%the blue and red crosses
%correspond to the minimum-time and minimum-energy trajectories,
%and the green cross corresponds to a compromise between these two.

%Note that as the time horizon increases past ,
%the energy consumed starts to increase.
%(This point is marked with a .)
%This is due to the fact that the engine characteristic
%has a positive ``idling loss'' $\gamma$.
%Because our model assumes the engine is always on,
%this idling loss continues over the entire horizon,
%even if the destination has been reached.

\paragraph{Sample trajectories\pp}
We now describe in detail the three trajectories
shown in figure~\ref{f-trajectories}.
We begin by describing the trajectory shown in green,
found by solving (\ref{e-opt-ctrl})
with a travel time of $T = 280$ seconds.
%This trajectory represents a trade-off between minimizing
%energy consumption and travel time.
The trajectory begins with acceleration at
the maximal rate $a_t^{\rm max}$;
the drive power is then decreased,
until it reaches zero, and the vehicle coasts
until the upper speed limit is reached
during time period $t=50$ seconds.
At this point the vehicle brakes,
then supplies constant power to maintain speed
at the speed limit.
After the speed limit ends ($t=100$ s),
the vehicle accelerates
to a roughly constant ``cruising speed'' of around $80$ $\rm km/h$.
Of particular interest is the existence of a coasting interval
during the last $50$ seconds,
during which the drive power is zero
as the vehicle coasts to the desired displacement.
This is a solution to the apparent dilemma
that a positive vehicle speed is required to reach the desired displacement,
yet any leftover positive kinetic energy at this point
can be considered wasted energy.
(One might argue that
this leftover kinetic energy could have been put to better use accelerating the vehicle earlier on;
evidently, this argument is false.)

The trajectory shown in blue is the minimum-time trajectory,
\ie, it is the smallest value of $T$ for which (\ref{e-opt-ctrl})
is feasible.
Note that this control depletes the internal energy exactly
as the desired displacement is reached.
The trajectory often accelerates aggressively,
using a substantial amount of power to get to a high speed quickly.
Even so, the minimum-time trajectory
has a coasting period that begins around $t =210$ seconds
and lasts until the end of the time period.

Finally, in red, we see the minimum-energy trajectory,
which is found by computing the value of $T$
that maximizes the optimal value of problem (\ref{e-opt-ctrl}).
As one might expect, the speed is kept lower than for the previous two trajectories,
which both decreases the amount of energy required to accelerate the vehicle,
and the power lost to drag.
However,
the minimum-energy trajectory still reaches the desired displacement
in a finite amount of time.
This is because our model assumes that the engine is turned on during the
interval $[0,T]$,
and idling incurs a power loss of $\gamma$.
Therefore, the minimum-energy control is motivated not to waste any time
in reaching the goal,
\ie, being fast also helps reduce wasted energy.

%We also note that
%none of the plotted values are linear,
%even over short intervals.
%(Except, of course,
%when the speed and acceleration limits are active.)
%This means that simple heuristics,
%like those mentioned in \cite{XXX},
%are likely to be suboptimal,
%even for the simple model considered here.

\section{Conclusions}
We used a simple optimization model to capture
the trade-off between vehicle energy consumption and travel time.
Several interesting extensions are possible,
especially for different drivetrain architectures,
and our formulation easily accommodates
drivetrain models formulated in terms of energy and power.
Modeling a solar car is a particularly simple extension,
which involves adding a time-varying prediction
of generated solar power
into the internal energy dynamics equation (\ref{e-energy-dyn}).
Another extension involves adding a time-dependent disturbance
to the kinetic energy dynamics (\ref{e-energy-dyn}),
which could model the predicted power loss (or gain) from traversing hilly terrain.

\section*{Acknowledgement}
The author would like to thank Stephen Boyd for useful input,
as well as Ashe Magalhaes and Gawan Fiore
for their insightful discussions of solar car racing.

\ifieee
\else
    \newpage
\fi

\bibliography{opt_racing}

\end{document}